\documentclass[12pt,leqno,a4paper]{article}
\usepackage{mathrsfs}

\usepackage{amsfonts}
\usepackage{amsmath, amssymb, amsthm, amsfonts}

\numberwithin{equation}{section}

\def\hf{{\textstyle{\frac12}}}
\def\a{\alpha}\def\b{\beta}

\def\e{\varepsilon}

\def\G{\Gamma} \def\g{\gamma}

\def\s{\sigma}
\def\z{\zeta}

\def\={\;=\;}

\def\zt{\zeta(\hf+it)}

\def\leq{\leqslant}
\def\geq{\geqslant}

\begin{document}
\baselineskip=17pt
\title{On a hybrid fourth moment involving the Riemann zeta-function}
\author{ Aleksandar Ivi\'c and Wenguang Zhai}
\date{}

\footnotetext[0]{Wenguang Zhai is supported by
 the National Key Basic Research Program of China (Grant No. 2013CB834201),
the National Natural Science Foundation of China (Grant No. 11171344),  the Natural
Science Foundation of Beijing (Grant No. 1112010) and the Fundamental Research Funds for the
Central Universities in China (2012Ys01).} \maketitle
\begin{center}
{\bf Abstract}
\end{center}
{\small
We provide explicit ranges for $\s$ for which the asymptotic formula
\begin{equation*}
\int_0^T|\zt|^4|\zeta(\sigma+it)|^{2j}dt \;\sim\; T\sum_{k=0}^4a_{k,j}(\sigma)\log^k T
\quad(j\in\mathbb N)
\end{equation*}
holds    as $T\rightarrow \infty$, when $1\leq j \leq 6$, where $\z(s)$ is
the Riemann zeta-function. The obtained ranges improve on
an earlier result of the authors. An application to a divisor problem
is also given.}

\bigskip

\section{\bf Introduction}

Let as usual $\zeta(s)=\sum\limits_{n=1}^\infty n^{-s}\ (\Re s > 1)$ denote
the Riemann zeta-function, where $s = \sigma+it$ is a complex
variable. Mean values of $\zeta(s)\,$ in the so-called  ``critical strip"
$\frac 12 \leq \sigma \leq 1$ represent a central topic in the theory of the
zeta-function (see e.g., the monographs \cite{IVi1} and \cite{IVi2}
   for an extensive account). Of special
interest are the moments on the so-called  ``critical" line $\sigma=\frac 12$.
Unfortunately as of yet no bound of
the form
\begin{equation}
\int_0^T |\zt|^{2m}dt \;\ll_{\varepsilon,m}\;T^{1+\varepsilon}\qquad (m\in {\mathbb N})
\end{equation}
is known to hold when $m\geq 3$, while in the cases $m = 1,2$ precise
asymptotic formulas for the integrals in question are known (see op. cit.).
Throughout this paper, $ \varepsilon$ denotes  fixed small positive    constants, not necessarily
the same ones at each occurrence, while $\ll_{a,\ldots}$ denotes the dependence of the
$\ll$-constant on $a,\ldots\,$.

\medskip
Having in mind the difficulties of establishing (1.2) when $m\geq3$, it appeared interesting
to consider the following problem.
For any fixed integer $j\geq 1,$ let $\sigma_{4,j}^{*}\;(\geq \hf)$ denote the infimum of all
$\sigma$ for which the estimate
\begin{equation}
\int_0^T|\zt|^4|\zeta(\sigma+it)|^{2j}dt\ll_{j,\varepsilon}T^{1+\varepsilon}
\end{equation}
holds. The left-hand side of (1.2) may be called a ``hybrid" moment, since it combines
moments on the lines $\Re s = \hf$ and $\Re s =\s$.
The problem is to estimate $\sigma_{4,j}^{*}$ for a given $j\in\mathbb N$.
If the well-known Lindel\"{o}f hypothesis ($\zt \ll_\e |t|^\e$) is true,
then $\sigma_{4,j}^{*}=\hf$  for any $j\geq 1.$ However, up to now even
$\sigma_{4,1}^{*}=\hf$ is out of reach by the use of existing  methods.
We cannot have $\sigma_{4,j}^{*}< \hf$ in view of the functional equation
$$
\z(s) \;=\; \chi(s)\z(1-s),\quad \chi(s) := \frac{\G(\hf(1-s))}{\G(\hf s)}\pi^{s-1/2}
\asymp |t|^{1/2-\s}.
$$

In his work \cite{IVi3} the first author investigated the integral in (1.2) for
the case $j=1$ and the case $j=2.$
 In particular, he proved that $\sigma_{4,1}^{*} \leq \frac{5}{6} = 0.8\bar3$,  while if $(k,\ell)$ is an exponent
pair (see e.g., \cite{GK} or Chapter 2 of \cite{IVi1} for definitions) with $3k+\ell<1$, then
$$
 \sigma_{4,2}^{*}\leq \max \left(\frac{\ell-k+1}{2},
\frac{11k+\ell+1}{ 8k+2}\right),
$$
which implies that $\sigma_{4,2}^{*}  \leq 1953/1984 =
0.984375$. Since $\z(\s+it) \ll \log |t|$ for $\s\geq1$, it is trivial
that $ \sigma_{4,j}^{*} \leq1$ for any fixed $j\in {\mathbb N}$.
At the end of   \cite{IVi3} it was stated, as an open problem,  to prove
the strict inequality  $ \sigma_{4,j}^{*} <1$   for any fixed $j\in {\mathbb N}$.

In  \cite{IZ}, which is a continuation of \cite{IVi3}, the authors proved that
indeed $ \sigma_{4,j}^{*} <1$ holds  for any fixed $j\in {\mathbb N}.$
  In fact,  if $(k,\ell)$ is an exponent pair with $\ell +(2j-1)k<1, $ then
we showed that $$\sigma_{4,j}^{*}\leq \frac{\ell+(6j-1)k}{1+4jk}.$$ In particular, we have
$\sigma_{4,2}^{*}\leq \frac{37}{38}=0.97368\cdots.$

In \cite{IZ} we also considered the possibilities of sharpening (1.2) to an asymptotic formula.
We showed that, for any given integer $j\geq1$, there exists a number $\sigma_1 =  \sigma_1(j)$
for which $\frac{3}{4} < \sigma_1 < 1$
such that, when $\sigma > \sigma_1$, there exists an asymptotic formula
for the integral in (1.2). This is
\begin{equation}
\int_0^T|\zt|^4|\zeta(\sigma+it)|^{2j}dt \;\sim T
\;\sum_{k=0}^4a_{k,j}(\sigma)\log^k T\quad(T\to\infty),
\end{equation}
where all the coefficients $a_{k,j}(\sigma) $, which depend on $\sigma$
and $j$, may be evaluated explicitly. However, in \cite{IZ} we did not provide explicitly
 the range of $\sigma$ for which (1.3) holds.

In this paper we shall provide some explicit values of $\s$ for which (1.3) holds.

{\bf Theorem 1.} {\it  The asymptotic formula} (1.3) {\it holds in the following ranges:
\begin{eqnarray*}
&&\sigma>\frac{4}{5}=0.8 \ \  (j=1),\\
&&\sigma> 0.904391\cdots\ \ ( j=2),\\
&&\sigma> 0.940001\cdots\ \ ( j=3),\\
&&\sigma> 0.959084\cdots\ \ ( j=4),\\
&&\sigma> 0.970734 \cdots\ \ ( j=5),\\
&&\sigma> 0.978286 \cdots\ \ ( j=6),\\
\end{eqnarray*}
 }
{\bf Corollary 1.}{\it We have
\begin{eqnarray*}
&&\sigma_{4,1}^{*}\leq \frac{4}{5}=0.8 \ \  (j=1),\\
&&\sigma_{4,2}^{*}\leq  0.904391\cdots\ \ ( j=2),\\
&&\sigma_{4,3}^{*}\leq  0.940001\cdots\ \ ( j=3),\\
&& \sigma_{4,4}^{*}\leq  0.959084\cdots\ \ ( j=4),\\
&& \sigma_{4,5}^{*}\leq  0.970734\cdots\ \ ( j=5),\\
&& \sigma_{4,6}^{*}\leq  0.978286\cdots\ \ ( j=6).\\
\end{eqnarray*} }

As an application of Theorem 1, we shall consider a weighted divisor problem.
Suppose that $\ell\geq 1$ is a fixed integer and $a$ is a fixed real number. Define
the divisor function
\begin{equation}
d_{4,\ell}(n) = d_{4,\ell}(n;a) =\sum_{n=n_1n_2}d_4(n_1)d_\ell(n_2)n_2^{-a},
\end{equation}
where $d_{k}(n)$ denotes the number of ways $n$ can be written as a product of
$k$ factors (so $d_k(n)$ is generated by $\z^k(s)$).
If $a=0,$ then $d_{4,\ell}(n)\equiv d_{4+\ell}(n).$

Henceforth we  consider only the case $ a>0.$ Suppose $X\geq2.$  It is expected that the summatory
function  $\sum_{n\leq X}d_{4,\ell}(n)$ is asymptotic to
$$
X\sum_{k=0}^{3}c_{k,\ell}(a)\log^k X+X^{1-a}\sum_{k=0}^{\ell-1}c_{k,\ell}^{\prime}(a)\log^k X
$$
as $X\to\infty,$ where the constants $c_{k,\ell}$ and $c_{k,\ell}^{'}$ are effectively
computable. More precisely, if one defines
$$
E_{4,\ell}(X):=\sum_{n\leq X}d_{4,\ell}(n)-X\sum_{k=0}^{3}c_{k,\ell}(a)\log^k X-X^{1-a}
\sum_{k=0}^{\ell-1}c_{k,\ell}^{\prime}(a)\log^k X,
$$
then we expect $E_{4,\ell}(X) = o(X)$ to hold as $X\to \infty$. Thus $E_{4,\ell}(X)$ should
represent the error term in the asymptotic formula for $\sum_{n\leq X}d_{4,\ell}(n)$.
It is also clear that the difficulty of the estimation of $E_{4,\ell}(x)$ increases with
$\ell$, and it also increases as $a$ in (1.4) gets smaller.

By using (1.2) and the  complex contour integration method, we can prove

\bigskip

{\bf Theorem 2.} {\it  If $\;\max\Bigl(\sigma^{*}_{4,j_0}-\hf, \hf-\frac{1}{\ell}\Bigr)\leq   a<\hf,$ then
for $\ell\geq1$ fixed we have
\begin{equation}
E_{4,\ell}(X) \;\ll_\e\; x^{1/2+\varepsilon},
\end{equation}
where $j_0=\hf\ell$ if $\ell$ is even, and $j_0=\hf(\ell+1)$ if $\ell$ is odd.}

\bigskip
From Theorem 2 and Corollary 1 we obtain at once

{\bf Corollary 2.} {\it  The estimate }(1.5) {\it holds for
\begin{eqnarray*}
&&   \textstyle\frac{3}{10}<a< \hf, \ \  (\ell=1,2),\\
&&   0.404391\cdots<a<\hf,\ \ ( \ell=3,4),\\
&&   0.440001\cdots<a<\hf,\ \ (\ell=5,6),\\
&&    0.459084\cdots<a<\hf,\ \ ( \ell=7,8),\\
&&   0.470734\cdots<a<\hf,\ \ ( \ell=9,10),\\
&&    0.478286\cdots<a<\hf,\ \ ( \ell=11,12).\\
\end{eqnarray*} }

\section{\bf The necessary  lemmas}

In order to prove our results, we require some  lemmas which will be given in this section.
The first lemma is the following
upper bound for the fourth moment of $\zeta(\frac 12+it),$ weighted by a Dirichlet polynomial.

\medskip
{\bf Lemma 2.1}. Let $a_1, a_2, \ldots, a_M$ be complex numbers. Then we have,
for $ \varepsilon > 0, M\geq 1$ and $T\geq 1$,
\begin{equation}
\int_0^T |\zt|^4\Bigl|\sum_{m\leq M}a_mm^{it}\Bigr|^2dt
\ll_\varepsilon T^{1+\varepsilon}M(1+M^2T^{-1/2})\max_{m\leq M}|a_m|^2.\end{equation}

\medskip This result is due to N. Watt \cite{Wa}. It is founded on the earlier works of
 J.-M. Deshouillers and H. Iwaniec \cite{DI}, which involved the use of Kloosterman sums, but
 Watt's result is sharper.

We also need some results on power moments of $\zeta(s).$

\medskip

{\bf Lemma 2.2}. For any fixed $A\geq 4,$ let us define $M(A)$ as
\begin{eqnarray*}
M(A)=\left\{\begin{array}{ll}
\frac{A-4}{8},&\mbox{if $4\leq A\leq 12,$}\\
\\
\frac{3A-14}{22},& \mbox{if $12\leq A\leq 178/13=13.6923\cdots,$}\\
\\
\frac{416A-2416}{2665},&\mbox{if $178/13\leq A\leq 20028/1313=15.253\cdots,$}\\
\\
\frac{7A-36}{48},&\mbox{if $20028/1313\leq A\leq 1836/101=18.178\cdots,$}\\
\\
\frac{32(A-6)}{205},&\mbox{if $A\geq 1836/101.$}
\end{array}\right.
\end{eqnarray*}
Then we have the estimate
\begin{equation}
\int_1^T |\zeta(\hf+it)|^Adt \;\ll_\e\; T^{1+M(A)+\varepsilon}.
\end{equation}

{\bf Proof.} The case $4\leq A\leq 178/13$ is contained in Theorem 8.2 of Ivi\'c \cite{IVi1}.
Now suppose that $A>178/13.$

Suppose that $t_1<t_2<\cdots<t_R$ are real numbers which satisfy
$$
|t_r|\leq T\ (r=1,2,\cdots, R),\ |t_s-t_r|\geq 1\quad (1\leq r\not= s\leq R),
$$
and
$$
|\zeta(\hf+it_r)|\geq V>0\qquad (r=1,2,\cdots, R).
$$

The large values estimate (8.29) of Ivi\'c \cite{IVi1} reads
 \begin{equation}
 \begin{split}
R&\ll TV^{-6}\log^8 T+T^{29/13}V^{-178/13}\log^{235/13} T\\
&\ll T^{29/13}V^{-178/13}\log^{235/13}T,
 \end{split}
 \end{equation}
if we note  that $ \zeta(\hf +it) \ll_\e |t|^{32/205+\varepsilon}$ (see M. N. Huxley
\cite{Hu1} and \cite{Hu2}).

We shall also use (8.33) of \cite{IVi1}, namely
 \begin{equation}
R\;\ll\; T^2V^{-12}\log^{16} T.
 \end{equation}

From (2.3) and (2.4) we obtain
\begin{eqnarray}
\ \ \ \ \ \;\;\int_1^T|\zeta(\hf+it)|^Adt\ll_\e\left\{\begin{array}{ll}
T^{2+\frac{3(A-12)}{22}+\varepsilon },&\mbox{if $12\leq A\leq 178/13,$}\\
T^{\frac{29}{13}+\frac{32(A-178/13)}{205}+\varepsilon },&\mbox{if $  A\geq 178/13.$}
\end{array}\right.
\end{eqnarray}

The formula (8.56) of Ivi\'c \cite{IVi1} reads
\begin{eqnarray}
R\ll\left\{\begin{array}{ll}
TV^{-6}\log^8 T ,&\mbox{if $V\geq T^{11/72}\log^{5/4}T,$}\\ \\
T^{15/4}V^{-24}\log^{61/2}T ,&\mbox{if $V< T^{11/72}\log^{5/4}T .$}
\end{array}\right.
\end{eqnarray}

From (2.4) and (2.6) we have
\begin{eqnarray}
\begin{split}
 & \int_1^T|\zeta(\hf+it)|^Adt\\
 &\ll_\e\left\{\begin{array}{ll}
T^{\max(1+\frac{32(A-12)}{205},2+\frac{7(A-12)}{48})+\varepsilon },&\mbox{if $12\leq A\leq 24,$}\\
T^{1+\frac{32(A-6)}{205}+\e },&\mbox{if $  A\geq 24.$}
\end{array}\right.
\end{split}
\end{eqnarray}

Now Lemma 2.2 for the case $A>178/13 $ follows from (2.5) and (2.7).

{\bf Lemma 2.3.}
For $1/2<\sigma <1$ fixed we define $m(\sigma)\;(\geq 4)$ as the
supremum  of all numbers $m \;(\geq 4)$ such that
\begin{equation}
 \int_1^{T}\left|\zeta(\sigma+it)\right|^{m}dt\ll_\e
T^{1+\varepsilon}
\end{equation}
 for any $\varepsilon>0.$ Then
$$\begin{array}{lll}
&m(\sigma)\geq 4/(3-4\sigma), & \frac{1}{2}< \sigma\leq \frac{5}{8},\\[8pt]
&m(\sigma)\geq 10/(5-6\sigma), &\frac{5}{8}\leq \sigma\leq \frac{35}{54},\\[8pt]
&m(\sigma)\geq 19/(6-6\sigma), &\frac{35}{54}\leq \sigma\leq \frac{41}{60},\\[8pt]
&m(\sigma)\geq 2112/(859-948\sigma), &\frac{41}{60}\leq \sigma\leq \frac{3}{4},\\[8pt]
&m(\sigma)\geq 12408/(4537-4890\sigma), &\frac{3}{4}\leq \sigma\leq \frac{5}{6},\\[8pt]
&m(\sigma)\geq 4324/(1031-1044\sigma), &\frac{5}{6}\leq \sigma\leq \frac{7}{8},\\[8pt]
&m(\sigma)\geq 98/(31-32\sigma), &\frac{7}{8}\leq \sigma\leq 0.91591\ldots,\\[8pt]
&m(\sigma)\geq (24\sigma-9)/(4\sigma-1)(1-\sigma), &
0.91591\ldots\leq \sigma\leq 1-\varepsilon .\end{array}$$

{\bf Proof.} This is Theorem 8.4 of Ivi\'c \cite{IVi1}. In Ivi\'c--Ouellet \cite{IVO}
some improvements have been obtained. Thus, it was shown there that $m(\s) \geq 258/(63-64\s)$
for $14/15 \leq \s \leq c_0$ and $m(\s) \geq (30\s-12)/(4\s-1)(1-\s)$ for $c_0 \leq \s \leq 1-\e$,
where $c_0 = (171+\sqrt{1602})/222 = 0.95056\cdots\,$.

{\bf Lemma 2.4.} Let $q\geq 1$ be an integer, $Q=2^q.$ Then for $|t|\geq 3$ we have
 \begin{equation}
 \zeta\left(1-\frac{q+2}{2^{q+2}-2}\right)\;\ll\; |t|^{1/(2^{q+2}-2)}\log |t|.
 \end{equation}
 We also have
 \begin{equation}
 \zeta(\textstyle\frac{5}{7}+it)\;\ll_\e\; |t|^{0.07077534\cdots+\varepsilon}\qquad(|t|\geq 2).
 \end{equation}

{\bf Proof.} The formula (2.9) is Theorem 2.12 of Graham and Kolesnik \cite{GK}. The estimate
(2.10) is to be found on page 66 of \cite{GK}. It improves (2.9) in the case when $q=2$, when
one obtains the exponent $\frac{1}{14} = 0.0714285\cdots\,$.

{\bf Lemma 2.5.} Suppose $\hf\leq \sigma_1<\sigma_2\leq 1$ are two real numbers such that
$$\zeta(\sigma_j+it)\;\ll_\e\; |t|^{c(\sigma_j)+\varepsilon}\qquad  (j=1,2),$$
then for $  \sigma_1\leq \sigma\leq \sigma_2  $ we have
\begin{equation}
\zeta(\sigma+it)\;\ll_\e\; {|t|}^{c(\sigma_1)\frac{\sigma_2-\sigma}{\sigma_2-\sigma_1}+c(\sigma_2)
\frac{\sigma-\sigma_1}{\sigma_2-\sigma_1}+\varepsilon}.
\end{equation}

{\bf Proof.} This follows from the well-known Phragm\'en-Lindel\"{o}f principle (convexity);
see e.g., Section 8.2 of \cite{IVi1}.

{\bf Lemma 2.6.} Let
$$
I(h,k) := \int_{-\infty}^\infty \left(\frac{h}{k}\right)^{-it}\zeta(\hf+\alpha+it)
\zeta(\hf+\beta+it)\zeta(\hf+\gamma-it)\zeta(\hf+\delta-it)w(t)dt,
$$
where $h,k\in\mathbb N, (h,k)=1$, and $\alpha, \beta, \gamma, \delta$ are
complex numbers $\ll 1/\log T$. Then for
$hk \leq T^{2/11-\e}$ we have
\begin{equation}
\begin{split}
I(h,k) &= \frac{1}{\sqrt{hk}}\int_{-\infty}^{\infty}w(t)\Biggl\{Z_{\a,\b,\g,\delta,h,k}(0)
+ \left(\frac{t}{2\pi}\right)^{-\a-\b-\g-\delta}Z_{-\g,-\delta,\a,-\b,h,k}(0)\\&
+ \left(\frac{t}{2\pi}\right)^{-\a-\g}Z_{-\g,\b,-\a,\delta,h,k}(0)
+ \left(\frac{t}{2\pi}\right)^{-\a-\delta}Z_{-\delta,\b,-\g,-\a,h,k}(0)
\\&
+ \left(\frac{t}{2\pi}\right)^{-\b-\g}Z_{\a,-\g,-\b,\delta,h,k}(0)
+ \left(\frac{t}{2\pi}\right)^{-\b-\delta}Z_{\a,\delta,\g,-\b,h,k}(0)\Biggr\}dt
\cr&
+ O_\e\Bigl(T^{3/4+\e}(hk)^{7/8}(T/T_0)^{9/4}\Bigr).
\end{split}
\end{equation}
The function $Z_{\ldots}(0)$ is given in term of explicit, albeit complicated
Euler products.

\bigskip

Formula (2.12) is due to   C. P. Hughes and M. P. Young \cite{HY}. It is intended
primarily for the asymptotic evaluation of the integral
\begin{equation}
\int_0^T|\zeta(\hf+it)|^4|M(\hf+it)|^2dt,
\end{equation}
where
$$
M(s) := \sum_{h\leq T^\theta}a(h)h^{-s}
$$
is a Dirichlet polynomial of length $T^\theta$ with coefficients $a(h)\;(\in \mathbb C)$.
The integral in (2.13) reduces to a sum of integrals of the type $I(h,k)$ after one
develops $|M(\hf+it)|^2 $ and chooses suitably the weight function $w(t)$, which is discussed
below. In general, the evaluation of the integral in (2.13)
is an important problem in  analytic number theory. It was studied by J.-M. Deshouillers
and H. Iwaniec \cite{DI}, N. Watt \cite{Wa}   and most recently by Y. Motohashi \cite{Mo},
all of whom used powerful methods from the spectral theory of the non-Euclidean Laplacian.
In \cite{HY}   Hughes and Young  obtained an asymptotic formula for (2.13) when
$\theta=\frac{1}{11}-\varepsilon.$
Two of the chief ingredients in their proof
  are an approximate functional equation for the product of
four zeta values, and the so-called ``delta method'' of Duke, Friedlander and Iwaniec
\cite{DFI}.  Watt's result (2.1) gives the expected upper bound $O_\e(T^{1+\e})$
in the range $\theta \leq \frac{1}{4}$, but does not produce an asymptotic formula
for the integral in (2.13) (or (2.1)). At the end of \cite{Mo}, Y. Motohashi comments on the value
$\theta=\frac{1}{11}-\varepsilon$ of \cite{HY}. He says:
``Our method should give a better result than theirs, if it is combined
with works by N. Watt on this mean value.''

\medskip

Note that the  bound $O_\e(T^{1+\e})$ for (2.13) with $T^\theta, \theta = \hf$
would give the hitherto unproved sixth moment of zeta-function in the form
$$
\int_0^T|\zt|^6dt \;\ll_\varepsilon\; T^{1+\varepsilon},
$$
which is (1.1) with $m=3$.

\medskip
The weight function $w(t)\;(\geq0)$ which appears in the integral in (2.12)
is a smooth function majorizing or minorizing the characteristic
function of the interval $[T, 2T]$. The fact that the integrand in (2.13) is
non-negative makes this effective. We shall actually take two such functions:
$w(t) = w_1(t)$ supported in $[T-T_0, 2T +T_0]$ such that $w_1(t) = 1$ for
$t\in [T, 2T]$, and $w(t) = w_2(t)$ supported in $[T, 2T]$ such
$w_2(t) = 1$ for $t\in [T-T_0, 2T-T_0]$.
For an explicit construction of such a smooth function $w(t)$ see e.g., Chapter 4
of the first author's monograph \cite{IVi2}. We then have, in either case,
$w^{(r)}(t) \ll_r T_0^{-r}$ for all $r = 0,1,2,\ldots$\,,
where $T_0$ is a parameter which satisfies $T^{1/2+\e} \ll T_0 \ll T$,
and appears in the error term in (2.12).

\section{\bf Proof of Theorem 1}

\subsection{The case when $j=1$}

\par
In this subsection we shall prove Theorem 1 in the case when $j=1.$
However, we shall deal with the general case and
restrict ourselves to $j=1$ only at the end of the proof.

Suppose $T\geq 10$. It suffices  to evaluate the integral
$$
\int_T^{2T}|\zt|^4|\zeta(\sigma+it)|^{2j}dt,
$$
replace then $T$ by $T2^{-j}$ for $j = 1,2,\ldots$ and sum the resulting
estimates. For convenience, henceforth we set
${\cal L}:=\log T.$
Let $s=\sigma+it,$ $\hf<\sigma\leq 1$ and $T\leq t\leq 2T.$
We begin with the well-known Mellin inversion integral (see e.g.,
the Appendix of \cite{IVi1}),
\begin{equation}
e^{-x}=\frac{1}{2\pi i}\int_{(c)}x^{-w}\Gamma(w)dw\qquad (c>0,\, x>0),
\end{equation}
where $\int_{(c)}$ denotes integration over the line $\Re w=c.$

Suppose $T^{1/11}\ll Y\ll T$ is a parameter to
 be determined later. In (3.1) we set $x=n/Y,$
multiply by $d_j(n)n^{-s}$ and then sum over $n.$ This gives
\begin{equation}
\sum_{n=1}^\infty d_j(n)e^{-n/Y}n^{-s}=\frac{1}{2\pi i}\int_{(2)}Y^w\zeta^j(s+w)\Gamma(w)dw.
\end{equation}

Suppose  $\s_0$ is fixed number which satisfies $\hf\leq \sigma_0<\min(1,\sigma)$
and will  be  determined later.
In (3.2) we shift the line of integration to $\Re w=\sigma_0-\sigma$ and apply the residue theorem.
The pole at $w=1-s$, which is  of degree $j$,  contributes the residue  which is $\ll T^{-10},$
by Stirling's formula for $\Gamma(w).$ The pole at $w=0$ contributes the residue
$\zeta^{j}(s).$ Thus we have
\begin{eqnarray}
&&\ \ \ \ \zeta^j(s)\\
&&=\sum_{n=1}^\infty d_j(n)e^{-n/Y}n^{-s}-\frac{1}{2\pi i}\int_{(\sigma_0-\sigma)}Y^w\zeta^j(s+w)\Gamma(w)dw+O(T^{-10})\nonumber.
\end{eqnarray}
By the well-known elementary estimate
$$\sum_{n\leq u}d_j(n)\ll u\log^{j-1}u$$
and partial summation it is easy to see that
$$\sum_{n> Y{\cal L}^2 } d_j(n)e^{-n/Y}n^{-s}\ll T^{-10}.$$
By Stirling's formula for $\Gamma(w)$ again we have
$$
\frac{1}{2\pi i}\int_{\Re w=\sigma_0-\sigma,|\Im w|>{\cal L}^2 }Y^w\zeta^j(s+w)\Gamma(w)dw\ll T^{-10}.
$$
Let $Y_1:=T^{1/11-\varepsilon}.$ Inserting the above two estimates into (3.3) we can write
\begin{equation}
\zeta^j(s)=B_1(s)+B_2(s)+B_3(s)+B_4(s),
\end{equation}
say, where
\begin{eqnarray}
&&B_1(s):=\sum_{n\leq Y_1}  d_j(n)e^{-n/Y}n^{-s},\nonumber\\
&&B_2(s):=\sum_{Y_1<n\leq Y{\cal L}^2}  d_j(n)e^{-n/Y}n^{-s},\nonumber\\
&&B_3(s):=-\frac{1}{2\pi i}\int_{\Re w=\sigma_0-\sigma,|\Im w|\leq{\cal L}^2 }
Y^w\zeta^j(s+w)\Gamma(w)dw,\nonumber\\
&&B_4(s):=O(T^{-10}).\nonumber
\end{eqnarray}
The partitioning in (3.4) is a new feature in the approach to this problem.
The flexibility is present in the parameters $Y$ and $\s_0$, which will allow us
to use Lemma 2.2 and Lemma 2.3, hence to connect our problem to the power moments
of $|\z(\s+it)|$.

\medskip

Therefore from (3.4) we have, since $|ab| \leq \hf|a|^2+\hf|b|^2$,
\begin{eqnarray*}
|\zeta(\sigma+it)|^{2j}&=&|B_1(\sigma+it)|^2\\
& +&\sum_{2\leq k\leq 4}O\left(|B_1(\sigma+it)B_k(\sigma+it)|+|B_k(\sigma+it)|^2\right).\nonumber
\end{eqnarray*}
Multiplying the above relation by $|\zt|^4$ and integrating, we obtain
\begin{eqnarray}
\ \ \ \ \ \ \ \ \ \int_T^{2T}|\zt|^4|\zeta(\sigma+it)|^{2j}dt=J_1+\sum_{2\leq k\leq 4}O(J_k+ J_k^{\prime}),
\end{eqnarray}
say, where
\begin{eqnarray*}
&&J_k: =\int_T^{2T}|\zt|^4|B_k(\sigma+it)|^2dt,\nonumber\\
&&J_k^{\prime}:= \int_T^{2T}|\zt|^4|B_1(\sigma+it)B_k(\sigma+it)|dt.
\end{eqnarray*}

The main contribution to the integral in (1.3) will come from the integral $J_1$,
with our choice $Y_1 = T^{1/11-\e}$.
In \cite{IZ}, the authors evaluated the integral  similar to
 $J_1$ with the help of the result of
  Hughes and  Young (Lemma 2.6). Actually, disregarding the harmless factor $e^{-n/Y}$,
 the integral $I_1$ in (4.6) of \cite{IZ}
 is just the integral $J_1$ if the
 parameter $Y=T^{1/(11j)-\varepsilon_1}$ therein is replaced by $Y_1=T^{1/11-\varepsilon}$ defined
 above. For the sake of completeness we shall give the details of the evaluation of $J_1$.
As a technical convenience,
we consider instead of $J_1$ the weighted integral
\begin{equation}
J^{*} := \int_{-\infty}^{\infty}w(t)|\z(\hf+it)|^4
\Bigl|\sum_{n\leq Y_1}d_j(n)e^{-n/Y}n^{-\s-it}\Bigr|^2\,dt,
\end{equation}
with $w(t) = w_j(t) \;(\geq 0; j=1,2)$ as in the discussion following Lemma 2.6.
We note that
\begin{equation*}
\int_{-\infty}^{\infty}w_2(t)|\z(\hf+it)|^4
\Bigl|\sum_{n\leq Y_1}\cdots\Bigr|^2\,dt \leq
J_1 \leq \int_{-\infty}^{\infty}w_1(t)|\z(\hf+it)|^4
\Bigl|\sum_{n\leq Y_1}\cdots\Bigr|^2\,dt,
\end{equation*}
and we shall show that the same asymptotic formula holds for the integral with $w_1(t)$ and
$w_2(t)$ above, which will show then that such a formula holds for  $J_1$ as well.
We write the square of the sum in (3.6) as
\begin{equation}
\begin{split}
&\Bigl|\sum_{n\leq Y_1}d_j(n)e^{-n/Y}n^{-\s-it}\Bigr|^2\\
= &\sum_{m,n\leq Y_1}d_j(m)d_j(n)e^{-m/Y}e^{-n/Y}\Bigl(\frac{m}{n}\Bigr)^{-it}(mn)^{-\s}
\\
= &\sum_{\delta\le Y_1}\delta^{-2\s}\sum_{h\leq Y_1/\delta,k\le Y_1/\delta,(h,k)=1}
d_j(\delta h)d_j(\delta k)e^{-\delta h/Y}e^{-\delta k/Y}(hk)^{-\s}\left(\frac{h}{k}\right)^{-it},
\end{split}
\end{equation}
where we put $m =\delta h, n = \delta k, (h,k)=1$.
With the aid of (3.7) it follows that $J^{*}$ reduces to the summation of
integrals of the type
$$
I^{*}(h,k) := \int_{-\infty}^{\infty}w(t)|\z(\hf+it)|^4\left(\frac{h}{k}\right)^{-it}d t
\qquad((h,k)=1).
$$
We continue now the proof of Theorem 1, and we  multiply (2.12) by
$$d_j(\delta h)d_j(\delta k)e^{-h/Y}e^{-k/Y}(hk)^{-\s}$$
and insert the resulting expression in (3.6). The error term in (2.12) makes a
contribution which will be, since $d_j(n) \ll_\e n^\e$,
\begin{equation*}
\begin{split}
&
\ll_\e \sum_{\delta\leq Y_1}\delta^{\e-2\s}
\sum_{h\leq Y_1/\delta, k\leq Y_1/\delta}T^{3/4+\e}(hk)^{7/8-\s}(T/T_0)^{9/4}
\\&
\ll_\e T^{3/4+\e}Y_1^{15/4-2\s}(T/T_0)^{9/4}.
\end{split}
\end{equation*}
Note that  $Y_1^{\frac{15}{4}-2\s} < Y_1^{\frac{11}{4}}$ because
$\s>\hf$. Therefore we see, since $Y_1 = T^{\frac{1}{11}-\e}$, as in
the discussion made in \cite{HY}, that we obtain first the desired asymptotic formula,
with an error term $O(T^{1-\e_1})$ for some $\e_1>0$, for the twisted integral $J^{*}$ in (3.6), with
$|\zt|^4$ replaced by
\[
\zeta(\hf+\alpha+it)
\zeta(\hf+\beta+it)\zeta(\hf+\gamma-it)\zeta(\hf+\delta-it).
\]
  Finally,
if $\alpha, \beta, \gamma, \delta$ all tend to zero, we obtain the desired
asymptotic  formula
 \begin{equation}
J_1\sim
 T\sum_{k=0}^4 b_{k;j}(\sigma)\log^k T\qquad(T\rightarrow\infty,),
 \end{equation}
 and the coefficients $b_{k;j}(\sigma)$ depend on $\sigma$ and $j.$
 It remains then to show that the contribution of $J_k$ and $J_k'$ in (3.5), for $2\leq k\leq4$,
 is of a lower order of magnitude than the right-hand side of (3.8), and Theorem 1 will follow.

\bigskip

We shall  estimate  the integral $J_2$ by Lemma 2.1.
We split the range of summation in $B_2(s)$ into $O(\log T)$
ranges of summation of the form
\[
Y_1 \leq M < n \leq M' \leq 2M \ll Y{\cal L}^2.
\]
 Hence by Lemma 2.1 and the well-known elementary bound
\begin{equation}
d_j(n)\;\ll_\e\; n^{\varepsilon}
\end{equation}
we have
\begin{equation*}
\begin{split}
  J_2&\ll {\cal L}\max_{Y_1\ll M \ll Y{\cal L}^2}\int_T^{2T}|\zt|^4
 \Bigl|\sum_{M<n\leq M'\leq 2M}  d_j(n)e^{-n/Y}n^{-\sigma-it} \Bigr|^{2}dt \\
&\ll_\varepsilon \max_{Y_1\ll M \ll Y{\cal L}^2}T^{1+\varepsilon}M(1+M^2T^{-1/2})
\max_{M<n\leq 2M}d_j^2(n)e^{-2n/Y}n^{-2\sigma}\nonumber\\
&\ll_\varepsilon  \max_{Y_1\ll M \ll Y{\cal L}^2}T^{1+\varepsilon}(M^{1-2\sigma}+M^{3-2\sigma}T^{-1/2})
 \\
&\ll_\e {\cal L}^2T^{1+3\varepsilon}Y_1^{1-2\sigma}+T^{1/2+3\varepsilon}Y^{3-2\sigma}{\cal L}^{8-4\sigma}
 \nonumber\\
&\ll_\e T^{1+\varepsilon}Y_1^{1-2\sigma}+T^{1/2+\varepsilon}Y^{3-2\sigma}. \\
\end{split}
\end{equation*}
Since $Y_1 = T^{\frac{1}{11}-\e}$, we see that
\begin{equation}
J_2\;\ll_\e\; T^{1-\varepsilon}
\end{equation}
if
\begin{equation}
Y \;=\; T^{\frac{1}{6-4\s}-\e},
\end{equation}
and the condition $T^{\frac{1}{11}} \ll Y \ll T$ is seen to hold.

We turn now to the estimation of the integral $J_3$. From its definition we have
$$
B_3(\sigma+it)\;\ll\; Y^{\sigma_0-\sigma}\int_{-{\cal L}^2 }^{{\cal L}^2 }|\zeta(\sigma_0+it+iv)|^j dv,
$$
hence by using this bound  and Cauchy's inequality we infer that
$$
|B_3(\sigma+it)|^2\;\ll\; Y^{2\sigma_0-2\sigma}{\cal L}^{2}
\int_{-{\cal L}^2 }^{{\cal L}^2 }|\zeta(\sigma_0+it+iv)|^{2j} dv.
$$
Thus by integration we have
$$
J_3
\ll Y^{2\sigma_0-2\sigma} {\cal L}^{2} \int_T^{2T}|\zt|^4
\left(\int_{-{\cal L}^2 }^{{\cal L}^2 }|\zeta(\sigma_0+it+iv)|^{2j} dv\right) dt.
$$

Suppose now that $\sigma_0$, besides $\hf \leqslant \s_0 < \min(1,\s)$,  also satisfies the condition
\begin{equation}
m(\sigma_0)\;>\;2j.
\end{equation}

Let
\begin{equation}
q: =\frac{m(\sigma_0)}{2j}, \ \ p:= \frac{m(\sigma_0)}{m(\sigma_0)-2j}.
\end{equation}
Then
$$p>1,\; q>1,\;  \frac{1}{p}+\frac{1}{q}=1,$$
and  by H\"{o}lder's inequality for integrals we obtain
\begin{equation}
\begin{split}
J_3&\ll Y^{2\sigma_0-2\sigma} {\cal L}^{2}\left(\int_T^{2T}|\zt|^{4p}dt\right)^{\frac 1p}\\
&\times
\left(\int_T^{2T}\left(\int_{-{\cal L}^2 }^{{\cal L}^2 }
|\zeta(\sigma_0+it+iv)|^{2j} dv\right)^q dt\right)^{\frac 1q}.
\end{split}
\end{equation}

We have
\begin{equation}
\int_T^{2T}|\zt|^{4p}dt  \;\ll_\e\; T^{1+M(4p)+\varepsilon},
\end{equation}
where we shall use the bounds for $M(A)$ furnished by Lemma 2.2.
By H\"{o}lder's inequality again we have
\begin{eqnarray*}
&&\left(\int_{-{\cal L}^2 }^{{\cal L}^2 }|\zeta(\sigma_0+it+iv)|^{2j} dv\right)^q\\
&&\ll \int_{-{\cal L}^2 }^{{\cal L}^2 }|\zeta(\sigma_0+it+iv)|^{2jq} dv \times
 \left(\int_{-{\cal L}^2 }^{{\cal L}^2 }1 dv\right)^{\frac qp}\\
 &&\ll {\cal L}^{\frac{2q}{p}}\int_{-{\cal L}^2 }^{{\cal L}^2 }|\zeta(\sigma_0+it+iv)|^{2jq} dv.
\end{eqnarray*}
Therefore
\begin{equation}
\begin{split}
&\int_T^{2T}\left(\int_{-{\cal L}^2 }^{{\cal L}^2 }|\zeta(\sigma_0+it+iv)|^{2j} dv\right)^q dt\\
&\ll {\cal L}^{\frac{2q}{p}}\int_T^{2T}\left(\int_{-{\cal L}^2 }^{{\cal L}^2 }
|\zeta(\sigma_0+it+iv)|^{2jq} dv\right)dt\\
&={\cal L}^{\frac{2q}{p}}\int_{-{\cal L}^2 }^{{\cal L}^2 }dv \int_T^{2T}|\zeta(\sigma_0+it+iv)|^{2jq} dt\\
&={\cal L}^{\frac{2q}{p}}\int_{-{\cal L}^2 }^{{\cal L}^2 }dv \int_T^{2T}
|\zeta(\sigma_0+it+iv)|^{m(\sigma_0)} dt\\
&\ll_\e T^{1+\varepsilon}.
\end{split}
\end{equation}

From (3.14)--(3.16) and (3.11) we obtain
\begin{equation}
\begin{split}
J_3&\ll Y^{2\sigma_0-2\sigma} {\cal L}^{2}\left( T^{1+M(4p)+\varepsilon}\right)^{\frac 1p}
\left( T^{1+\varepsilon}{\cal L}^{\frac{2q}{p}+2} \right)^{\frac 1q}\\
&\ll Y^{2\sigma_0-2\sigma}T^{1+\frac{M(4p)}{p}+\varepsilon}{\cal L}^4\ll_\e T^{1-\e}
\end{split}
\end{equation}
if $\;T^{ \frac{M(4p)}{p}+3\varepsilon} \ll Y^{2\s-2\s_0}$. With the choice (3.11)
this condition reduces to
\begin{equation}
\s \;>\: \frac{\frac{3M(4p)}{p} + \s_0}{\frac{2M(4p)}{p}+1}.
\end{equation}

\medskip

To bound the integrals $J_k^{\prime}$ (see (3.5)) note that
from (3.8), (3.10), (3.17) and Cauchy's inequality for integrals we obtain
\begin{equation}
J_k^{\prime}\leq J_1^{1/2}J_k^{1/2}\ll T^{1-2\varepsilon}{\cal L}^2\ll T^{1-\varepsilon}
\qquad  (k=2, 3).
\end{equation}

Obviously we have
\begin{equation}
J_4\ll T^{-18},
\end{equation}
and  consequently
\begin{equation}
J_4^{\prime}\ll T^{-16}.
\end{equation}

From (3.8), (3.10), (3.17) and (3.19)--(3.21) we obtain that, if (3.12) and (3.18) hold,
\begin{equation*}
 \int_T^{2T}|\zt|^4|\zeta(\sigma+it)|^{2}dt\;\sim\;
 T\sum_{k=0}^4 b_{k;j}(\sigma)\log^k T\qquad(T\to\infty),
 \end{equation*}
This implies that
\begin{equation*}
 \int_1^{T}|\zt|^4|\zeta(\sigma+it)|^{2}dt\;\sim\;
 T\sum_{k=0}^4 a_{k;j}(\sigma)\log^k T\qquad(T\to\infty),
 \end{equation*}
where the $a_{k;j}$'s are constants which are easily expressible in term of the $b_{k;j}$'s.

Now we determine the permissible range of $\sigma $ from (3.18) for the case $j=1.$
We take $\sigma_0=\frac{5}{8}$. Lemma 2.3 gives $m(\s_0) = m(\frac{5}{8}) \geq8$,
so (3.11) holds, and $p = \frac{4}{3}$. Then  (3.18) reduces to
$\sigma>\frac{4}{5}.$

\medskip
{\bf Remark 1.}
When $j=2, 3, 4,$ the above procedure can also give non-trivial results.
Actually, when $j=2,$ we take $\sigma_0=\frac{35}{54}$ and (3.18) becomes
$\sigma>\frac{71}{78}=0.91025\cdots.$
 When $j=3,$ we take $\sigma_0=\frac{5}{6}$ and (3.18) becomes
$\sigma>\frac{659}{690}=0.95507\cdots.$
When $j=4,$ we take $\sigma_0=\frac{7}{8}$ and (3.18) becomes
$\sigma>\frac{221}{229}=0.96506\cdots.$ However, in Subsection 2 we shall give better ranges
for $\s$ in these three cases.

{\bf Remark 2.}
When $j > 4$, the above method does not give good results in view of the
existing bounds for the functions $M(A)$ and $m(\s)$ defined in Lemma 2.2 and Lemma 2.3 respectively. However, in that
case it is not difficult to see that (1.3) holds for $\s > \s^{*}_{4,j}$,
the infimum of numbers for which (1.2) holds. Thus (1.3) will hold for
$\s > (\ell + (6j-1)k)/(1+4jk)$ when $(k,\ell)$ is an exponent pair.
To see this, note first that the discussion preceding (3.12) yields
\begin{equation}
J_3 \ll Y^{2\sigma_0-2\sigma} {\cal L}^{2} \max_{|v|\leq {\cal L}^{2}}
\int_T^{2T}|\zt|^4|\zeta(\sigma_0+it+iv)|^{2j} \,dt.
\end{equation}
This is the almost the same integral as the initial one, and the
conclusion of Theorem 1 of our joint paper \cite{IZ} holds, namely
\begin{equation}
\int_0^T|\zt|^4|\z(\s+it+iv)|^{2j}\,dt \ll_{j,\e} T^{1+\e}\quad(|v| \leq {\cal L}^2),
\end{equation}
if
\begin{equation}
\s > \s_0 \= \frac{\ell+(6j-1)k}{1+4jk},\qquad
\ell +(2j-1)k<1
\end{equation}
and $(k,\ell)$ is an exponent pair. With $\s_0$ as in (3.23) and $\s \geq \s_0+\delta\;(\delta>0)$
one has trivially $J_3 \ll T^{1-1/11}$ for $\delta,\e$ sufficiently small
(since $Y^{2\s_0-\s} \ll T^{-2/11+2\delta\e}$), and we
get an asymptotic formula for the initial integral in the range
$\s > \s_0$ for $j>4.$

\medskip

{\bf Remark 3.}
We may further discuss the asymptotic formula (1.3).
Denote by, say, $E(T;\sigma,j)$ the
difference between the left and right-hand side in (1.3), thus
$E(T;\sigma,j)$ is the error term in the asymptotic
formula for our integral. Let $c(\sigma,j)$ be the infimum
of numbers $c$ such that, for a given $j\in \mathbb N$,
$$
E(T;\sigma,j) \ll T^c.
$$
We know that $c(\sigma,j) < 1$ by \cite {IZ}, and it seems reasonable to
expect that $c(\sigma,j) \geq \hf$. Namely in case when $j=0$, we have the fourth
moment of $|\zt|$, and in this case a precise asymptotic formula is known, and the
exponent of the error term cannot be smaller than $\hf$ (see \cite{IViM}).
However, obtaining any qualitative results on $c(\sigma,j)$ will be difficult,
one of the reasons being that it is hard from the method of Hughes and Young \cite{HY}
to get explicit $O$-estimates for the error terms in their formulas.

\subsection{The case when $j\geq 2$}\

To deal with  the case $j\geq 2$ we shall use an induction method.
Namely, for each $j\geq 1,$ we shall prove that there is a constant $\hf<c_j<1$ such that
(3.23) holds in the range $\sigma>c_j.$ When $j=1,$ we can take 
$c_1=\frac{4}{5}$ from the result in Subsection 1.

Let $C(\sigma)>0$ be a function which connects $
(\frac{5}{7},0.07077534\cdots)$ and the   points $(a_q,b_q)\ (q\geq 3)$
with line segments, where
\begin{equation}
a_q:=1-\frac{q+2}{2^{q+2}-2}, \ b_q: = \frac{1}{2^{q+2}-2},
\end{equation}
and $q = q(j)$ will be suitably chosen.
 We then have, in view of Lemma 2.4 and Lemma 2.5,
\begin{equation}
 \zeta(\sigma+it)\;\ll_\e\; |t|^{C(\sigma)+\varepsilon}\quad (\sigma\geq \textstyle\frac{5}{7}).
 \end{equation}

 Now we suppose that $j\geq 2 $
 and   we have already defined  $c_{l}$ for any $1\leq l<j.$ From (3.22) and (3.26) we have
\begin{equation}
\begin{split}
J_3&\ll Y^{2\sigma_0-2\sigma} {\cal L}^{2} \max_{|v|\leq {\cal L}^{2}}
\int_T^{2T}|\zt|^4|\zeta(\sigma_0+it+iv)|^{2+2(j-1)} \,dt\\
& \ll_\e Y^{2\sigma_0-2\sigma}T^{2C(\sigma_0)+2\varepsilon} {\cal L}^{2} \max_{|v|\leq {\cal L}^{2}}
\int_T^{2T}|\zt|^4|\zeta(\sigma_0+it+iv)|^{ 2(j-1)} \,dt\\
&\ll_\e Y^{2\sigma_0-2\sigma}T^{2C(\sigma_0)+1+3\varepsilon }  \\
&\ll_\e T^{\left(\frac{1}{6-4\sigma}-\varepsilon\right)(2\sigma_0-2\sigma)+2
C(\sigma_0)+1+3\varepsilon }\\
&\ll_\e T^{ \frac{  \sigma_0- \sigma}{3-2\sigma}   +2C(\sigma_0)+1+3\varepsilon }
\end{split}
\end{equation}
if $\sigma_0>c_{j-1}.$

Take $\sigma_0=c_{j-1}+\delta,$ where $\delta>0$ is a small positive constant.
When
\begin{equation}
\sigma>\frac{6C(c_{j-1})+c_{j-1}}{4C(c_{j-1})+1},
\end{equation}
  from (3.27) we have
$$J_3\;\ll_\e\; T^{1-\varepsilon}$$
if $\delta, \varepsilon$ are sufficiently small.

We define the sequence $c_j\;(j\geq 1)$ as follows:
\begin{equation}
c_1=\frac{4}{5},\ \ \  c_j:= \frac{6C(c_{j-1})+c_{j-1}}{4C(c_{j-1})+1}
\qquad (j\geq 2).
\end{equation}
It is easy to see  that $c_j<1$ for any $j$ since $C(\sigma)<\hf(1-\sigma).$  From the above procedure and the
results in Subsection 1 we see that the asymptotic formula (1.3) holds  for $\sigma>c_j$ for any
$j\geq 2.$

We provide now the explicit values of $c_j $ when $j=2, 3, 4, 5, 6,$ and we remark that continuing in
this fashion we could obtain the values for $j>6$ as well.

1. The case $j=2$:  from Lemma 2.4 we have $C(\frac{5}{7})=0.07077534\cdots,$
$C(\frac{5}{6})=\frac{1}{30}=0.03333333\cdots.$ Thus from Lemma 2.5 we have
$C(\frac{4}{5})=0.0438170952\cdots.$ Hence
$$c_2=\frac{6C(4/5)+4/5}{4C(4/5)+1}=0.904391\cdots.$$

\medskip

2. The case $j=3$: from (3.25) we have $a_4=\frac{28}{31}<c_2<a_5=\frac{119}{126}.$ From Lemma 2.4
 we have $C(\frac{28}{31})=\frac{1}{62},$ $C(\frac{119}{126})=\frac{1}{126}.$ From Lemma 2.5 we get
$C(c_2)=0.01589736\cdots.$ Hence
\begin{equation*}
\begin{split}
c_3&= \frac{6C(c_2)+c_2}{4C(c_2)+1}
 = 0.9400013\cdots
\end{split}
\end{equation*}

\medskip
3. The case $j=4$: we have $a_4=\frac{28}{31}<c_3<a_5=\frac{119}{126}.$   From Lemma 2.5 we get
$C(c_3)=0.008819601\cdots.$ Hence
\begin{equation*}
\begin{split}
c_4&= \frac{6C(c_3)+c_3}{4C(c_3)+1}
=0.959084\cdots.
\end{split}
\end{equation*}

\medskip
4. The case $j=5$: we have $a_5=\frac{119}{126}<c_4<a_6=\frac{123}{127}.$   From Lemma 2.5 we get
$C(c_4)=0.005502913\cdots.$ Hence
\begin{equation*}
\begin{split}
c_5&= \frac{6C(c_4)+c_4}{4C(c_4)+1}
= 0.970734\cdots.
\end{split}
\end{equation*}

\medskip
5. The case $j=6$: we have $a_6=\frac{123}{127}<c_5<a_7=\frac{501}{510}.$   From Lemma 2.5 we get
$C(c_5)= 0.0035902\cdots.$ Hence
\begin{equation*}
\begin{split}
c_6&= \frac{6C(c_5)+c_5}{4C(c_5)+1} = 0.978286\cdots.
\end{split}
\end{equation*}

{\bf Remark 4.} It is not difficult to evaluate additional values of $c_j$. For example, we have
$c_7=0.983536,$ $ c_8=0.987254,$ $ c_9=0.990005,$ $ c_{10}=0.992046,$ $ c_{11}=0.993616.$
When $j$ large, the value of $c_j$ is close to $1.$

{\bf Remark 5.} The values of $c_j\;(j\geq 2)$ depend on the upper bound   of
$\zeta(\sigma+it).$ Therefore  we can improve the values of $c_j\;(j\geq 2)$ 
if we have better upper bounds   for $\zeta(\sigma+it).$  
For example,  instead of Lemma 2.4 (Th. 2.12 of Graham-Kolesnik \cite{GK}), we could use Theorem 4.2
of theirs (p. 38), which is strong for any $q\geq 1.$ Then we can get small improvements for any
$j\geq 2.$ We also remark that we have (see (7.57) of \cite{IVi1})
$$
\z(\s+it) \;\ll\; t^{(k+\ell-\s)/2}\log t\qquad(\s \geq \hf, \ell - k \geq\s),
$$
where $(k,\ell)$ is an exponent pair. A judicious choice of the exponent pair $(k,\ell)$,
especially the use of new exponent pairs due to M.N. Huxley (see e.g., his papers \cite{Hu2} and \cite{Hu3}),
would likely lead to some further small improvements. 

  Kevin Ford \cite{F} proved
$$|\zeta(\sigma+it)| \leq 76.2t^{4.45(1-\sigma)^{3/2}}\log^{2/3}t$$
for $\hf \leq\sigma\leq 1, t \geq 3$.
This estimate is quite explicit, and best when $\sigma$ is close to $1.$
This estimate would imply better values of $c_j$ when $j$ is large. There is, however, no 
simple procedure which yields (in closed form) the range for $\s$ for which the asymptotic
formula (1.3) holds, for any given $j$.

\bigskip
\section{\bf Proof of Theorem 2}

In this section we shall prove Theorem 2. By the definition of  the generalized function
$d_k(n)$ we have, for $0\leq a < \hf$ and $\Re s>1$,
\begin{equation}
\sum_{n=1}^\infty d_{4,\ell}(n)n^{-s}= \sum_{n_1,n_2=1}^\infty d_4(n_1)d_\ell(n_2)n_2^{-a}(n_1n_2)^{-s}
=\zeta^4(s)\zeta^\ell (s+a).
\end{equation}

By using Perron's inversion formula (see e.g., the Appendix of \cite{IVi1}) we have
\begin{equation}
\sum_{n\leq X}d_{4,\ell}(n)=\frac{1}{2\pi i}\int_{1+\varepsilon-iX}^{1+\varepsilon+iX}\zeta^4(s)
\zeta^\ell (s+a)\frac{X^s}{s}ds+O_\e(X^\varepsilon)
\end{equation}
if we note that $d_{4,\ell}(n)\ll_{\e,\ell} n^\varepsilon.$ Now we put
$j_0=\hf\ell$ if $\ell$ is even, and $j_0=\hf(\ell+1)$ if $\ell$ is odd,
and then move the line of integration in (4.2) to $\sigma=\hf.$
In doing this we encounter two poles. These are $s=1$, a pole of order four, and $s=1-a$ which is a pole of order
$\ell.$ It is easy to verify that the sum of residues of the integrand in (4.2) is of the form
(1.4). Thus from (1.5), (4.2) and   the residue theorem we obtain
\begin{equation}
E_{4,\ell}(X)=I_1+I_2-I_3+O_\e(X^\varepsilon),
\end{equation}
say, where
\begin{eqnarray*}
&&I_1:=\frac{1}{2\pi i}\int_{\frac{1}{2}-iX}^{\frac{1}{2}+iX}\zeta^4(s)\zeta^\ell (s+a)\frac{X^s}{s}ds,\\
&&I_2:=\frac{1}{2\pi i}\int_{\frac{1}{2}+ iX}^{1+\varepsilon+iX}\zeta^4(s)\zeta^\ell (s+a)\frac{X^s}{s}ds,\\
&&I_3:=\frac{1}{2\pi i}\int_{\frac{1}{2}-iX}^{1+\varepsilon-iX}\zeta^4(s)\zeta^\ell (s+a)\frac{X^s}{s}ds.
\end{eqnarray*}

For $\zeta(s)$ we have the bounds
\begin{eqnarray}
\zeta(\sigma+it)\ll\left\{\begin{array}{ll}
(2+|t|)^{\frac{1-\sigma}{3}}\log (2+|t|),&\mbox{if $\hf\leq \sigma\leq 1,$}\\ \\
\log (2+|t|) ,& \mbox{if $1\leq \sigma\leq 2.$}
\end{array}\right.
\end{eqnarray}
which follows from the standard bounds
$$
\zt \ll |t|^{\frac{1}{6}}\log|t|, \quad\z(\s+it) \ll \log|t|\quad(\s\geq 1, |t|\geq2)
$$
and convexity (see e.g., (1.67) of \cite{IVi1}).
Recalling the condition
$$
\max\Bigl(\sigma^{*}_{4,j_0}-\frac{1}{2}, \frac{1}{2}-\frac{1}{\ell}\Bigr)\leq   a<\frac{1}{2},
$$
we obtain
\begin{equation}
\begin{split}
& |I_2|+|I_3|\\&\ll \int_{\frac{1}{2}}^{1-a}X^{\frac{4(1-\sigma)}{3}+
\frac{\ell(1-\sigma-a)}{3}-1}X^\sigma\log^{4+\ell} Xd\sigma\\
& +\int_{1-a}^{1}X^{\frac{4(1-\sigma)}{3}-1}X^\sigma\log^{4+\ell} Xd\sigma+
\int_{1}^{1+\varepsilon}X^{-1}X^\sigma\log^{4+\ell} Xd\sigma\\
&\ll X^{\frac{1}{2}}\log^{4+\ell} X.
\end{split}
\end{equation}

Now we estimate $I_1.$ When $\ell$ is even,
we obtain directly $I_1 \ll_e X^{1/2+\e}$, since $j_0 = \hf\ell$.
Therefore we consider in detail the case when  $\ell$ is odd.
Let $j_1=\ell-j_0=\hf(\ell-1)=j_0-1.$ Then we have by Cauchy's inequality that
\begin{equation*}
\begin{split}
& \int_0^T|\zeta(\hf+it)|^4|\zeta(\hf+a+it)|^\ell dt\\
&=\int_0^T|\zeta(\hf+it)|^4|\zeta(\hf+a+it)|^{j_0+j_1} dt\\
&\ll \left(\int_0^T|\zeta(\hf+it)|^4|\zeta(\hf+a+it)|^{2j_0} dt\right)^{1/2} \\
& \times
 \left(\int_0^T|\zeta(\hf+it)|^4|\zeta(\hf+a+it)|^{2j_1} dt\right)^{1/2}.
 \end{split}
\end{equation*}
Since $\sigma^{*}_{4,j_0}-\hf \leq   a$, we have $\hf+a\geq \sigma^{*}_{4,j_0}.$ Hence from (1.2) we obtain
$$\int_0^T|\zeta( \hf+it)|^4|\zeta(\hf+a+it)|^{2j_0} dt\ll_\e T^{1+\varepsilon}.$$
Similarly we have
$$\int_0^T|\zeta(\hf+it)|^4|\zeta(\hf+a+it)|^{2j_1} dt\ll_\e T^{1+\varepsilon}.$$

From the above three estimates we obtain
\begin{eqnarray*}
  \int_0^T|\zeta(\hf+it)|^4|\zeta(\hf+a+it)|^\ell dt \ll_\e T^{1+\varepsilon},
\end{eqnarray*}
which combined with  integration by parts gives
\begin{equation}
I_1\;\ll_\e\; X^{\frac{1}{2}+\varepsilon}.
\end{equation}

\noindent
By combining (4.3), (4.5) and (4.6) we complete  the proof of Theorem 2.


\vfill
\eject

\vfill\eject
\bigskip

\bigskip\bigskip\bigskip

{\small
\noindent
 Aleksandar Ivi\'{c},\\
Katedra Matematike RGF-a Universiteta u Beogradu,\\
 Dju\v sina 7, 11000 Beograd,\\
 Serbia\\
E-mail: {\tt ivic@rgf.bg.ac.rs}\\

\noindent
Wenguang Zhai,\\
Department of Mathematics,
\\
China University of Mining and Technology,
\\
Beijing  100083, P.R.China\\
E-mail: {\tt zhaiwg@hotmail.com}
}
\end{document}